\documentclass[12pt,reqno,a4paper]{amsart}
\usepackage{amsmath,amssymb,amsfonts,amsthm,amscd}
\usepackage[mathscr]{eucal}
\usepackage{comment}

\usepackage{diagrams}  

\usepackage{hyperref}  

\date{13 (26) October 2015}

\author{Theodore Voronov}
\address{{School of Mathematics,  University of Manchester,    Manchester,   M13 9PL,  UK}
\newline
{\hspace{\parindent}}$\hphantom{m\;}$Department of Quantum Field Theory, Tomsk State University, Tomsk, 634050, Russia}
\email{theodore.voronov@manchester.ac.uk}

\title[On volumes of classical supermanifolds]{On volumes of classical supermanifolds}
\dedicatory{Dedicated to K.}

\newtheorem{theorem}{Theorem}
\newtheorem*{theor}{Theorem}
\newtheorem{prop}{Proposition}[section]
\newtheorem{cor}{Corollary}[section]
\newtheorem*{coro}{Corollary}
\theoremstyle{definition}
\newtheorem{example}{Example}[section]

\newtheorem*{Rem}{Remark}
\newtheorem*{Rems}{Remarks}

\def\co{\colon\thinspace}

\newcommand{\cal}{\EuScript}

\newcommand{\V}{\cal{V}}

\DeclareMathOperator{\Ber}{Ber}
\DeclareMathOperator{\Ker}{Ker}

\DeclareMathOperator{\Mat}{\mathbf{Mat}}

\renewcommand{\Im}{\mathop{\mathrm{Im}}}
\renewcommand{\Re}{\mathop{\mathrm{Re}}}

 \DeclareMathOperator{\diag}{diag}

 \DeclareMathOperator{\vol}{vol}
 \DeclareMathOperator{\sign}{sign}

\DeclareMathOperator{\ind}{ind}
\DeclareMathOperator{\indc}{ind_{\CC}}
\DeclareMathOperator{\dimc}{dim_{\CC}}
\DeclareMathOperator{\Arg}{Arg}

\newcommand{\der}[2]{{\frac{\partial {#1}}{\partial {#2}}}}

\newcommand{\D}[2]{{\frac{D{#1}}{D{#2}}}}

\newcommand{\lder}[2]{{\partial {#1}/\partial {#2}}}

\newcommand{\R}[1]{{\mathbf R}^{#1}}
\newcommand{\C}[1]{{\mathbf C}^{#1}}

\newcommand{\CC}{\mathbf C}

\newcommand{\ZZ}{{\mathbf Z}}

\newcommand{\p}{\partial}

\newcommand{\widebar}{\overline}

\renewcommand{\a}{\alpha}
\renewcommand{\b}{\beta}
\newcommand{\e}{\varepsilon}
\newcommand{\z}{\zeta}

\renewcommand{\O}{\Omega}
\renewcommand{\o}{\omega}

\newcommand{\G}{{\Gamma}}
\newcommand{\h}{\eta}
\renewcommand{\t}{\theta}

\newcommand{\x}{{\xi}}
\renewcommand{\d}{\delta}
\renewcommand{\t}{\theta}

\renewcommand{\L}{{\Lambda}}

\newcommand{\ft}{{\tilde f}}

\newcommand{\at}{{\tilde a}}
\newcommand{\bt}{{\tilde b}}

\newcommand{\Xt}{{\tilde X}}

\newcommand{\Tt}{{\tilde T}}

\newcommand{\vp}{{\boldsymbol{v}}}

\newcommand{\up}{{\boldsymbol{u}}}
\newcommand{\xp}{{\boldsymbol{x}}}

\newcommand{\wed}{\wedge}

\newcommand{\sph}{\mathbf S}
\newcommand{\pro}{\mathbf{CP}}
\newcommand{\sti}{\mathbf{CV}}
\newcommand{\st}{\mathbf{V}}
\newcommand{\un}{\mathbf{U}}
\newcommand{\gra}{\mathbf{CG}}
\newcommand{\gr}{\mathbf{G}}

\begin{document}
\begin{abstract}
We consider the volumes of classical supermanifolds such as the supersphere, 
complex projective superspace, 
and Stiefel and Grassmann supermanifolds, 
with respect to the natural metrics or symplectic structures. We show that the formulas for the volumes, upon certain universal normalization, can be obtained by an analytic continuation from the  formulas for the volumes of the corresponding ordinary manifolds.

Volumes of nontrivial supermanifolds may identically vanish.
In 1970s,  Berezin discovered that the total Haar measure of the unitary supergroup $\un(n|m)$ vanishes unless $m=0$ or $n=0$, i.e., unless it reduces to the ordinary unitary group $\un(n)$ or $\un(m)$. Witten recently suggested that the (Liouville) volume of a compact even symplectic supermanifold should always be zero if it is not an ordinary manifold. Our calculations provide   counterexamples  to this conjecture. On the other hand, we give a simple explanation of Berezin's statement  and generalize it to the Stiefel supermanifold $\st_{r|s}(\C{n|m})$. There are also  possible connections with the recent works  by Mkrtchyan and Veselov on `universal formulas' in Lie algebra theory.
\end{abstract}

\maketitle

\tableofcontents

\section{Introduction}

The concept of volume in the super world may show
unexpected features. For example, the total volume of a supermanifold
may vanish.  This is due to the fact that the underlying ``Berezin
integral'' is not of the usual measure-theoretic nature, so
considerations based on positivity of measure are not applicable. A striking
example is the superanalog of the unitary group, the unitary supergroup
$\un(n|m)$, whose total ``Haar measure" is zero whenever $nm>0$ (Berezin,
1970s). (This reflects also in the properties of the representation theory of the supergroup $\un(n|m)$. One cannot divide by zero volume as well as by an infinite volume. It has been noticed that the presence of odd dimensions is somewhat similar to noncompactness, see, e.g.~\cite{tv:gitnew}.)

In a discussion in fall 2012,  E.~Witten asked me  whether the Liouville volume
of every compact symplectic supermanifold should be zero if it is not an ordinary manifold. Although  he suggested  serious arguments in support of this conjecture,  in a short time I   provided   a counterexample. As such one can take the
superanalog of a complex projective space endowed with the analog of the
classical Fubini--Study form. The simplest case when here the Liouville supervolume is nonzero  is $\pro^{1|1}$, for which $\vol(\pro^{1|1})=2\pi$ (and there is a remarkable scaling invariance, see details in the main text).\footnote{\,In fact, as I later discovered by a thorough search of the literature, the answer for the volume of $\pro^{n|m}$ had  been previously known~\cite{lazaroiu:berezintoeplitz},  probably not so widely, otherwise Witten's question would not have been raised.}

The investigation   of the volume of the complex projective superspace has led me to observations that are interesting in their own right. Namely, the explicit formula for the volume of
the complex projective superspace $\pro^{n|m}$ is an analytic
continuation of the corresponding formula for the ordinary complex
projective space $\pro^n$ up to a universal normalization factor that depends only on  dimension. This triggered a study of other examples. As it turns out, this is also the case   for
other classical supermanifolds  such as the superspheres, the complex Stiefel and Grassmann supermanifolds. Formulas for the volumes of these supermanifolds can be obtained by simple geometric considerations.
In particular, in the course of this study we have obtained a simple explanation of the Berezin statement concerning the supergroup $\un(n|m)$ together with a generalization to the Stiefel supermanifold (for which the volume does not have to vanish identically).

The fact that volumes of classical manifolds and supermanifolds can be expressed by analytic functions, and very forms of these functions, is very remarkable. Here is a meeting point with the recent studies of ``universal formulas" in Lie algebra theory
by Mkrtchyan and Veselov (see, in particular,~\cite{mkrtchyanveselov:dualandnegat}, \cite{mkrtchyan:nonpertchernsimons}, \cite{mkrtchyanveselov:volumes}, \cite{mkrtchyan:quadrbarnes}). (One source of these ideas is~\cite{mkrtchyan:nandminusn}, another is the works of Vogel and Deligne.) There remain many open questions here.

The  analytic functions expressing the normalized volumes depend on complex variables whose geometric meaning are  ``indices'' (the numbers $n-m$ for superdimensions $n|m$). This  suggests that there may be a nontrivial extension of the theory to the infinite-dimensional case (and  some elliptic operators given) for which these volumes should still make sense.

The structure of the paper is as follows. In Section~\ref{sec.cpnm}, we consider the volumes of the supersphere and the complex projective superspace. In Section~\ref{sec.others}, we consider the complex Stiefel supermanifold,  the unitary supergroup and the Grassmann supermanifold. In Section~\ref{sec.further}, we briefly discuss the results. In Appendix~\ref{sec.prel}, we recall  some background information concerning supermanifolds, Berezin integration, and Riemannian submersions.

The starting point of this investigation was a question raised by E.~Witten, for which I am very grateful to him. The results of the work were discussed with many colleagues and reported  at seminars in Manchester, Sheffield and the Steklov Institute in Moscow, as well at the conferences in Bia{\l}owie\.{z}a (July 2013 and June 2014), Brussels (April 2014) and Skolkovo (October 2014). It is a pleasure to thank in particular V.~Buchstaber,  M.~Cahen,  R.~L.~Fernandes, S.~Gutt,  H.~Khudaverdian, K.~Mackenzie, R.~Mkrtchyan, Yu.~Neretin, A.~Odesskii,  A.~Odzijewicz, S.~Paycha,  A.~G.~Sergeev, A.~Veselov, and I.~Volovich.

\section{Superspheres and complex projective superspaces} \label{sec.cpnm}

\subsection{The supersphere and its volume}
To define the supersphere $\sph^{n|2m}$, consider the superspace $\R{n+1|2m}$  as a `vector supermanifold'  equipped with the   scalar product such that
\begin{equation*}
    (\xp,\xp)= (x^1)^2+\ldots+(x^{n+1})^2 + 2\xi^1\h^1+\ldots+2\xi^m\h^m\,,
\end{equation*}
in standard  coordinates $x^1,\ldots, x^{n+1},\x^1,\h^1,\ldots,\x^m,\h^m$.  (We need an even number of odd variables to ensure the nondegeneracy of the scalar product.) Then, as usual, the \emph{supersphere}  of radius $R$ is the submanifold $\sph^{n|2m}_R\subset \R{n+1|2m}$ of all   vectors of length $R$ and is given by
\begin{equation}\label{eq.sphere}
    (x^1)^2+\ldots+(x^{n+1})^2 + 2\xi^1\h^1+\ldots+2\xi^m\h^m= R^2\,.
\end{equation}
When $R=1$, we drop the reference to the radius.

The Euclidean metric on $\R{n+1|2m}$,
\begin{equation*}
    \d s^2 = (\d x^1)^2+\ldots+(\d x^{n+1})^2 + 2\d\xi^1\d\h^1+\ldots+2\d\xi^m\d\h^m\,,
\end{equation*}
induces a Riemannian metric on the submanifold $\sph^{n|2m}_R$, which gives rise to a volume element. When we speak about the volume of the supersphere, we refer to the total volume w.r.t. this Riemannian volume element.

\begin{theorem} The volume of the supersphere of radius $R$ equals
\begin{equation}\label{eq.volsphR}
    \vol(\sph^{n|2m}_R)= 2R^{n-2m}\,\frac{\pi^{\frac{n+1}{2}}\,2^m}{\G(\frac{n+1}{2}-m)}\,.
\end{equation}
\end{theorem}

\smallskip

{\small
Formula~\eqref{eq.volsphR} can be obtained from~\eqref{eq.sphere} by a calculation with the delta-function. Or, alternatively, it can be obtained via the Gaussian integral  (as in the ordinary case).

}

\smallskip
The expression for the volume of the supersphere has probably been established many times by different people; in particular, it was deduced among other things in~\cite{haunch:2007} by Adam Haunch (who was a joint Ph.~D. student of H.~M.~Khudaverdian and myself). When $m=0$, we have the familiar formula for the ordinary sphere. As in the ordinary case, it makes sense to look at the `even-dimensional' and `odd-dimensional' superspheres separately.\footnote{\,`Even-dimensional' and `odd-dimensional' refer here to the even part of the supermanifold dimension only. The pun is unavoidable, but unintended.} We have the following corollary (with a slight  change of   notation).
\begin{coro}
\begin{equation}\label{eq.volevsphR}
    \vol(\sph^{2n|2m}_R)=
    2R^{2n-2m}\,\frac{\pi^{n+\frac{1}{2}}\,2^m}{\G(n-m+\frac{1}{2})}\,.
\end{equation}
\begin{equation}\label{eq.volodsphR}
    \vol(\sph^{2n+1|2m}_R)
    = 2R^{2n+1-2m}\,\frac{\pi^{n+1}\,2^m}{\G(n-m+1)}\,.
\end{equation}
\end{coro}


Recall that the $\G$-function has poles at all negative integers and only at them.  From here we observe that  the volume of the `odd-dimensional' supersphere $\sph^{2n+1|2m}_R$   vanishes when $m>n$, and when it is nonzero, it can be written as
\begin{equation}\label{eq.volodsphRfact}
    \vol(\sph^{2n+1|2m}_R)
    = 2R^{2n+1-2m}\,\frac{\pi^{n+1}\,2^m}{(n-m)!}\,.
\end{equation}
At the same time, the volume of the `even-dimensional' supersphere $\sph^{2n|2m}_R$ never vanishes (and can always be expressed via double factorials).

\subsection{Recollection of $\pro^{n|m}$\,. Metric, symplectic structure  and volume element}
Let us recall the superanalogs of constructions familiar for ordinary complex projective space.
The complex projective superspace $\pro^{n|m}$ is defined exactly as in the ordinary case, as the quotient of the space of all nonzero (even) vectors in $\C{n+1|m}$ by the action of $\CC^*=\CC \setminus \{0\}$.

Let $z^a, \z^{\mu}$  be standard coordinates on $\C{n+1|m}$. They are also the `homogeneous coordinates' on $\pro^{n|m}= \C{n+1|m}\setminus \{0\}/ \CC^*$. 
For every fixed $k$ 
we can consider the variables $w^a_{(k)}:=\frac{z^a}{z^k}$ and $\t^{\mu}_{(k)}:=\frac{\z^{\mu}}{z^k}$ as  coordinates on the $k$th affine chart of $\pro^{n|m}$. They are also referred to as `inhomogeneous coordinates'. (Altogether there are $n+1$   charts. Henceforth the number $k$ is suppressed.)

We may   identify $\C{n+1|m}$ with $\R{2n+2|2m}$ by writing $z^a=x^a+y^a$ and $\z^{\mu}=\x^{\mu}+\h^{\mu}$ with real variables $x^a$, $y^a$, $\x^{\mu}$, and $\h^{\mu}$.
The standard Euclidean metric on $\C{n+1|m}$,
\begin{equation*}
    \d s^2 = \sum \d z^a\d \bar z^a + i \sum \d \z^{\mu} \d \bar \z^{\mu} = \sum \left((\d x^a)^2 +(\d y^a)^2\right) + 2\sum \d \x^{\mu} \d \h^{\mu}\,,
\end{equation*}
is the real part of the Hermitian form
\begin{equation*}
    H = \sum \d z^a\otimes \d \bar z^a + i \sum \d \z^{\mu}\otimes \d \bar \z^{\mu}\,.
\end{equation*}
The imaginary  part (up to a factor of $2$) is the  standard  symplectic form on $\C{n+1|m}$,
\begin{equation*}
    \o_0 =\frac{i}{2} \left(\sum dz^a d\bar z^a -i \sum d\z^{\mu}d\bar\z^{\mu}\right)=\sum dx^a dy^a + \frac{1}{2} \sum \left((d\x^{\mu})^2+ (d\h^{\mu})^2\right)\,.
\end{equation*}
Note that
\begin{equation*}
    \o_0 =\frac{i}{2}\, \p\bar\p N_0\,,
\end{equation*}
where
\begin{equation*}
    N_0=\sum z^a\bar z^a +i \sum \z^{\mu}\bar \z^{\mu}\,
\end{equation*}
is the scalar square of the radius-vector.

The symplectic form $\o_0$ does not descend on $\pro^{n|m}$. However, its restriction to the supersphere does. We have the analog of the familiar diagram
\begin{diagram}
\sph^{2n+1|2m}_R & \rInto^{\quad i\quad}  &   \C{n+1|m}\setminus \{0\} \\
& \rdTo_{p} &    \dTo_{q}    \\
& & \pro^{n|m} &
\end{diagram}
Here the slanted down-arrow is the analog of the Hopf fibration. The equation of the supersphere $\sph^{2n+1|2m}_R$ in complex coordinates is
\begin{equation*}
    \sum z^a\bar z^a +i \sum \z^{\mu}\bar \z^{\mu}=R^2\,.
\end{equation*}
The following proposition gives the  analog of the Fubini--Study form for $\pro^{n|m}$.

\begin{prop} The form $i^*\o_0$ on $\sph^{2n+1|2m}_R$ is $\widehat{U(1)}=\Pi T U(1)$-invariant and hence there exists a unique form $\o\in \O^2(\pro^{n|m})$   such that $i^*\o_0=p^*\o$.
\end{prop}

In local coordinates, the symplectic form $\o$ on $\pro^{n|m}$ is given by
\begin{equation*}
    \o=\frac{i}{2}\,R^2\frac{(dw\cdot d\bar w -id\t\cdot d\bar\t)N-(dw\cdot \bar w+id\t\cdot\t)(d\bar w\cdot w-id\bar \t\cdot \t)}{N^2}\,,
\end{equation*}
where $N=w\cdot \bar w + i \t\cdot \bar\t$. Here we sometimes write   sums such as $\sum dw^a d\bar w^a$, etc.,  in the form  $dw\cdot d\bar w=\sum dw^a d\bar w^a$. Note also that it is convenient to keep one `fake' inhomogeneous coordinate which is identically $1$, so that $N$ includes the summand $1$. The expression for $\o$ in homogeneous coordinates, which is the same as the pullback $q^*\o$ on $\C{n+1|m}\setminus \{0\}$, has similar form:
\begin{equation*}
    q^*\o=\frac{i}{2}R^2\frac{(dz\cdot d\bar z -id\z\cdot d\bar\z)N_0-(dz\cdot \bar z+id\z\cdot\z)(d\bar z\cdot z-id\bar \z\cdot \z)}{N_0^2}
\end{equation*}
where $N_0=z\cdot \bar z + i \z\cdot \bar\z$, as before. The form $\o$ can also be expressed via its K\"ahler potential:
\begin{equation*}
    \o=\frac{i}{2}\,R^2\p\bar\p \ln N\,, \quad q^*\o=\frac{i}{2}\,R^2\p\bar\p \ln N_0\,.
\end{equation*}
We loosely refer to both the symplectic form $\o$ and the corresponding Riemannian metric as to the `Fubini--Study form' on $\pro^{n|m}$. So everything is completely similar to the ordinary case. (The usual definition of the Fubini--Study form does not include any $R$, but for our purposes it is convenient to have a scale factor. When we wish to emphasize it in the formulas, we include $R$ as a subscript.)

The following volume element on $\pro^{n|m}$ is at the same time the Riemannian volume element and the Liouville volume element for the symplectic structure.

\begin{prop}
In local coordinates $w^a, w^a, \t^{\mu}, \bar\t^{\mu}$\,, 
\begin{equation}\label{eq.volelemcpnm}
    dV_{\pro^{n|m}_R} = R^{2(n-m)} \left(\frac{i}{2}\right)^{n-m}\frac{\left[dw^1, d\bar w^1,\ldots\,dw^n, d\bar w^n\,|\, d\t^1,d\bar\t^1,\ldots,d\t^m,d\bar\t^m\right]}{N^{n-m+1}}
\end{equation}
\end{prop}
Our next task is to calculate the total volume of  $\pro^{n|m}$ w.r.t. the volume element~\eqref{eq.volelemcpnm}. In order to do so, we establish the relation between the volume elements of $\pro^{n|m}_R$ and the supersphere.

\subsection{Volume elements for $\sph^{2n+1|2m}_R$ and  $\pro^{n|m}_R$}

To the find the desired relation, we     introduce local coordinates on $\C{n+1|m}$, with the same number of charts as for the projective superspace, which consist of the polar radius $r$, the angular coordinate $\a$ in the fibers of the Hopf fibration and the inhomogeneous coordinates on $\pro^{n|m}$. For each fixed $k$ (the chart number), we have
\begin{equation*}
    z^a=r e^{i\a} N^{-1/2}w^a\,, \quad \z^{\mu}=re^{i\a}N^{-1/2}\t^{\mu}\,.
\end{equation*}
Here $\a=\Arg z^k$, $r^2=z\cdot \bar z +i \z\cdot \bar\z$, $w^a=\frac{z^a}{z^k}$, and  $N=w\cdot \bar w + i \t\cdot \bar\t$ is as before.

\begin{prop}
In the coordinates $r,\a,w^{a},\bar w^{a}, \t^{\mu}, \bar \t^{\mu}$, the Euclidean volume element for $\C{n+1|m}$ is
\begin{equation*}
    dV_{\C{n+1|m}} = \left(\frac{i}{2}\right)^{n-m} r^{2n+1-2m} \frac{\left[dr, d\a, dw^1, d\bar w^1,\ldots\,dw^n, d\bar w^n\,|\, d\t^1,d\bar\t^1,\ldots,d\t^m,d\bar\t^m\right]}{N^{n-m+1}}\,.
\end{equation*}
\end{prop}

\begin{coro} For the supersphere $\sph^{2n+1|2m}_R$, in the local coordinates $\a,w^{a},\bar w^{a}, \t^{\mu}, \bar \t^{\mu}$,
\begin{align*}
    dV_{\sph^{2n+1|2m}_R} &= R^{2n+1-2m} \left(\frac{i}{2}\right)^{n-m}  \frac{\left[d\a, dw^1, d\bar w^1,\ldots\,dw^n, d\bar w^n\,|\, d\t^1,d\bar\t^1,\ldots,d\t^m,d\bar\t^m\right]}{N^{n-m+1}}\\
    &= Rd\a\cdot dV_{\pro^{n|m}_R}\,.
\end{align*}
\end{coro}
That means that we have a factorization of the volume element  in the Hopf fibration, hence the total volume  of the supersphere $\sph^{2n+1|2m}_R$ is the product of the volume of the complex projective superspace $\pro^{n|m}_R$ and the circumference of the circle $\sph^1_R$.

\subsection{Formula for the volume of $\pro^{n|m}$}

\begin{theorem}
 The volume of the complex projective superspace is given by the expression:
 \begin{equation}
  \vol(\pro^{n|m}_R)=R^{2(n-m)}\frac{\pi^n2^m}{\G(n-m+1)}\,.
 \end{equation}
\end{theorem}
\begin{example}
 $\vol(\pro^{1|1}_R)=2\pi$ (no dependence on radius). In particular, this is nonzero, which gives the simplest example of a compact symplectic supermanifold with a nonzero volume.
\end{example}

\begin{example}
 $\vol(\pro^{n|m}_R)=0$ if $m>n$ (because the $\G$-function has poles at all negative integers).
\end{example}

\subsection{A useful reformulation: volumes as analytic functions}

We can give the following useful reformulation of the obtained formulas for the volumes.
Recall that dimensions of supermanifolds take values in the ring $\hat\ZZ=\ZZ[\Pi]/(\Pi^2-1)$, so that $n|m=n+m\Pi$. Define
the   \emph{index} $\ind M^{n|m}$ of a supermanifold $M^{n|m}$ by
\begin{equation*}
    \ind M^{n|m}:=n-m\,,
\end{equation*}
if $\dim M^{n|m}=n|m$\,. We may say that $\ind M =\chi(\dim M)$, for a   ring homomorphism $\chi\co \hat\ZZ\to \ZZ$, $n|m \mapsto n-m$. For a complex supermanifold $M^{2n|2m}$ we can also speak about the dimension and index over $\CC$,
\begin{equation*}
    \ind_{\CC} M^{2n|2m}:=n-m\,,
\end{equation*}
where $\dim_{\CC}M^{2n|2m}=n|m$.
We also define  the \emph{Gaussian factor} $g_D$ for a given dimension $D=n|m \in \hat\ZZ$, which is the number
\begin{equation*}
    g_{n|m}:=(\sqrt{\pi})^n(\sqrt{2})^m\,.
\end{equation*}
(For $D=n|2m$, this is the value of the standard Gaussian integral over $\R{n|2m}$, see the Appendix. For aesthetical purposes, though we do not need that, we define $g_D$ for arbitrary dimensions $D$.)

Define for a complex variable $z$ the following functions:
\begin{align}
 \V(\sph; R,z)&:= R^z\,\frac{2\sqrt{\pi}}{\G(\frac{z+1}{2})}\\
 \V(\pro; R,z)&:=R^{2z}\,\frac{1}{\G(z+1)}\,.
\end{align}
(Here $R$ is a real parameter included for convenience.)
Then the formulas for the volumes of the superspheres and the complex projective superspaces can be re-written as:
\begin{equation}
 \vol(\sph^{n|2m}_R)=g_D\cdot \V(\sph; R,z)
\end{equation}
where $D=\dim \sph^{n|2m}=n|2m$, $z=\ind \sph^{n|2m}=n-2m$, and
\begin{equation}
 \vol(\pro^{n|m}_R)=g_{2D}\cdot \V(\pro; R,z)
\end{equation}
where $D=\dimc \pro^{n|m}=n|m$ (so the real dimension is $2D$), $z=\indc \pro^{n|m}=n-m$\,. It follows that if the volumes of the supermanifolds are re-scaled by the Gaussian factors for the corresponding dimensions, then upon such a normalization, they  are expressed by analytic functions of   ``index-type'' variables:
\begin{align}
 \frac{\vol(\sph^{D}_R)}{g_D}= \V(\sph; R,z) \ \text{\ where  $z=\ind \sph^D$\,,} \ \intertext{and}
 \frac{\vol(\pro^{D}_R)}{g_{2D}}= \V(\pro; R,z) \ \text{\ where  $z=\indc\pro^D$}\,.
\end{align}
The analytic functions $\V(\sph; R,z)$ and $\V(\pro; R,z)$ of the variable $z$ can then be regarded as  the  volumes of the spheres and complex projective spaces ``normalized'' in the above sense and analytically continued to arbitrary complex values of the  indices.

Since for ordinary manifolds   index coincides with   dimension, \emph{the normalized volumes of the supermanifolds in question   can be viewed  as the result of an analytic continuation of the normalized volumes of their classical counterparts}: first the dimension such as $n$ is replaced by a complex variable $z$ and then $z$ is substituted by $n-m$.

\begin{Rem}
The relation between the volume of the complex projective (super)space, the volume of the odd-dimensional (super)sphere and the circumference of the circle remains true for arbitrary complex $z$\,:
\begin{equation}
 \V(\sph; R,2z+1)= \V(\pro; R,z) \cdot \V(\sph; R,1)\,.
\end{equation}
\end{Rem}

\section{The unitary supergroup, and  Stiefel and Grassmann supermanifolds} \label{sec.others}

\subsection{The unitary supergroup}

Recall that the unitary supergroup $\un(n|m)$ is defined by the matrix equation
\begin{equation*}
    gHg^*=H
\end{equation*}
for a matrix $g\in\Mat(n|\,m\,; \CC)$,
where
\begin{equation*}
    g^*=\begin{pmatrix}
          g_{00} & g_{01} \\
          g_{10} & g_{11} \\
        \end{pmatrix}^*=
    \begin{pmatrix}
          \bar g_{00}^T & \bar g_{10}^T \vspace{3pt}\\
         - \bar g_{01}^T & \bar g_{11}^T
        \end{pmatrix}
\end{equation*}
and
\begin{equation*}
   H= \begin{pmatrix}
       1 & 0 \\
       0 & i
     \end{pmatrix}\,.
\end{equation*}
More geometrically, the supergroup $\un(n|m)$ can be defined as the supergroup of all linear operators on $\C{n|m}$ preserving the Hermitian metric. In complete analogy with the ordinary case, one can see that
\begin{align*}
    \dim \un(n|m)=(n|m)^2\,,\\
    \ind \un(n|m)=(n-m)^2\,.
\end{align*}
Note also that $\un(n|m)\cong \un(m|n)$\,.

The following remarkable fact was discovered by Berezin in early 1970s.
\begin{theor}[Berezin]
For $\un(n|m)$,  the total Haar measure is zero if both $n>0$ and $m>0$ \emph{(i.e., if it does not reduce to the ordinary unitary group)}.
\end{theor}

\begin{example} For the simplest case of $\un(1|1)$, the theorem can be checked directly. One can find the following explicit parametrization of the matrices $g=g(\a,\b\,|\,\t)\in \un(1|1)$\,:
\begin{equation*}
    g(\underbrace{\a,\b}_{\text{even real}} |\underbrace{\t}_{\text{odd complex}})= \begin{pmatrix}
                   e^{i\a}(1+\frac{i}{2}\t\bar\t) & \t \vspace{5pt}\\
                   i\bar\t e^{i\b} & e^{-i\a}(1-\frac{i}{2}\t\bar\t)e^{i\b}
                 \end{pmatrix}\,.
\end{equation*}
By a direct calculation,
\begin{equation*}
    [dg\cdot g^{-1}] = - 2i [d\a,d\b | d\t,d\bar\t]\,,
\end{equation*}
relative to the basis
\begin{equation*}
    e_1=\begin{pmatrix}
          i & 0 \\
          0 & 0 \\
        \end{pmatrix}\,, \quad
    e_2=\begin{pmatrix}
          0 & 0 \\
          0 & i \\
        \end{pmatrix}\,, \quad
    \e_1=\begin{pmatrix}
          0 & 1 \\
          i & 0 \\
        \end{pmatrix}\,, \quad
    \e_2=\begin{pmatrix}
          0 & i \\
          1 & 0 \\
        \end{pmatrix}\,
\end{equation*}
in $\mathfrak{u}(1|1)$. This gives (up to a constant factor) the invariant volume element on $\un(1|1)$. Therefore,
\begin{equation*}
    \vol \un(1|1) =0\,,
\end{equation*}
by the definition of the Berezin integral.
\end{example}

We shall calculate the volume of the complex Stiefel  supermanifold, for which the unitary supergroup is a particular case;  we shall obtain a simple geometric explanation of Berezin's theorem.

\subsection{The Stiefel  supermanifold}
Consider $\C{n|m}$ with the standard Hermitian form. For the complex Stiefel supermanifold, we use the notations $ \st_{r|s}(\C{n|m}) = \sti(n|m, r|s)$\,. It is defined as the closed submanifold
\begin{equation*}
    \st_{r|s}(\C{n|m}) \subset \underbrace{\C{n|m}\times \dots \times \C{n|m}}_{r}\times \underbrace{\Pi \C{n|m}\times \times \dots \times \Pi\C{n|m}}_{s}
\end{equation*}
of all orthonormal $r|s$-frames in $\C{n|m}$. (In greater detail, the vectors in the frame are orthogonal, the first $r$ vectors are even and of square $1$ and the last $s$ vectors are odd of square $i$.)
The Riemannian metric on $\st_{r|s}(\C{n|m})$ is induced from the Euclidean  metric on $\C{n|m}$ (the real part of the Hermitian form).

As in the ordinary case, $\st_{r|s}(\C{n|m})=\un(n|m)/\un(n-r|m-s)$. One can also find that
\begin{align*}
    \dim \st_{r|s}(\C{n|m})&= (r|s)\cdot \bigl(2(n|m)-(r|s)\bigr)= r(2n-r)+s(2m-s)\,|\, 2(ns+mr-rs)\,,\\
    \ind \st_{r|s}(\C{n|m}) &= (r-s) \bigl(2(n-m)-r+s\bigr)
\end{align*}
(generalizing the familiar formula $\dim \st_r(\C{n})=r(2n-r)$).

Note that if $V$ is a Hermitian vector superspace, then on the superspace $\Pi V$ with the revered parity is also naturally induced a Hermitian structure, by
\begin{equation*}
     \langle \Pi \up, \Pi \vp\rangle := i(-1)^{\tilde \up}\langle\up,\vp\rangle\,.
\end{equation*}
If the initial structure was positive-definite, the induced structure is also positive-definite.
From here follows a  useful statement.

\begin{prop} The   parity reversion functor $\Pi$ induces an isometry
\begin{equation*}
    \st_{r|s}(\C{n|m})\approx \st_{s|r}(\C{m|n})\,.
\end{equation*}
\end{prop}

The following statement is crucial for calculation of the volume of the Stiefel supermanifold.

\begin{theorem}
For $r>0$, there is a fiber bundle
\begin{equation*}
    \begin{CD}
    \st_{r-1|s}(\C{n-1|m}) @>>> \st_{r|s}(\C{n|m})\\
    @. @VVV \\
    @.   \sph^{2n-1|2m}
    \end{CD}\,.
\end{equation*}
(so that $\st_{r-1|s}(\C{n-1|m})$ is the standard fiber).
This fibration is a Riemannian submersion.
\end{theorem}

Recall that a fibration of Riemannian (super)manifolds is a Riemannian submersion    if the projection induces the  isometry of every horizontal subspace  onto the tangent space to the base. In such a case, there is a factorization of the volume element of the total space and in a homogeneous situation the volume of the total space is the product of the volumes of the base and the standard fiber  (see the Appendix).

The fiber bundle $\st_{r|s}(\C{n|m})\to \sph^{2n-1|2m}$ is not trivial, but on a dense domain it can be seen as a direct product decomposition and by iteration we may speak about a ``diffeomorphism'' (on a dense domain),  for $r>0$,
\begin{multline*}
    \st_{r|s}(\C{n|m})\approx \sph^{2n-1|2m} \times \sph^{2n-3|2m}\times \dots \times \sph^{2(n-r)+1|2m}
    \times \st_{0|s}(\C{n-r|m})\\
    \approx \sph^{2n-1|2m} \times \sph^{2n-3|2m}\times \dots \times \sph^{2(n-r)+1|2m}
    \times \st_{s}(\C{m|n-r})\approx \\
     \sph^{2n-1|2m} \times \sph^{2n-3|2m}\times \dots \times \sph^{2(n-r)+1|2m}
    \times \sph^{2(m-1)+1|2(n-r)} \times   \dots \times \sph^{2(m-s)+1|2(n-r)}
\end{multline*}
(the last equality for $s>0$). This is just for better visualization. In any case, we have the corresponding factorization for the volumes. Now, recall that the volume of the odd-dimensional  supersphere $\sph^{2k+1|2m}$ vanishes if $k-m<0$; therefore, the   volume
\begin{equation*}
    \vol(\st_{r|s}(\C{n|m}))= \vol(\sph^{2n-1|2m}) \cdot \vol(\sph^{2n-3|2m})\cdot \ldots \cdot \vol(\sph^{2(n-r)+1|2m})
    \cdot \vol(\st_{s}(\C{m|n-r}))
\end{equation*}
can be nonzero only if $n-r-m+1>0$, i.e., $r< n-m+1$. At the same time, if also  $s>0$, then
\begin{equation*}
    \vol(\st_{s}(\C{m|n-r}))=\vol(\sph^{2(m-1)+1|2(n-r)})\cdot     \ldots \cdot \vol(\sph^{2(m-s)+1|2(n-r)})\,,
\end{equation*}
which for the same reason can be nonzero only if $m-1-n+r+1>0$, i.e., $r> n-m$. This is impossible for an integer $r$. Therefore we arrive at the following theorem, the first part of which generalizes Berezin's statement about the volume of $\un(n|m)$.

\begin{theorem} 1. The volume of the complex Stiefel supermanifold  $\st_{r|s}(\C{n|m})$ vanishes unless $r=0$ or $s=0$.

2. For $s=0$, the volume is given by
  \begin{multline}
     \vol (\st_{r}(\C{n|m}))= \vol (\sph^{2(n-1)+1|2m}) \cdot \ldots \cdot  \vol (\sph^{2(n-r)+1|2m}) \\
     =g_D\,R^{\chi}\,(2\sqrt{\pi})^r\,\frac{1}{\G(n-m)\G(n-m-1)\dots \G(n-m-r+1)} 
  \end{multline}
  (where we have re-introduced the scale factor $R$ in the metric). Here
  $D=\dim \st_{r}(\C{n|m})$ and $\chi= \chi(D)=\ind \st_{r}(\C{n|m})$.
\end{theorem}

(We   used the multiplicativity of the Gaussian factor $g_{\dim}$ and of the scale factor $R^{\ind}$.)

By re-writing the fraction as
\begin{equation*}
    \frac{\G(n-m-r)\G(n-m-r-1)\dots \G(1)}{\G(n-m)\G(n-m-1)\dots \G(1)}\,,
\end{equation*}
we obtain the following corollary.

\begin{coro}[``Final formula'']
The volume of the complex Stiefel supermanifold when it is nonzero is given by
\begin{equation}
    \vol \st_{r}(\C{n|m})=  g_D \cdot  R^{\chi(D)} (2\sqrt{\pi})^r \frac{G(n-m-r+1)}{G(n-m+1)}\,,
\end{equation}
where $D=\dim \st_{r}(\C{n|m})= r(2n-r)|mr$, and $\chi(D)=\ind \st_{r}(\C{n|m})= r(2(n-m)-r)$\,.
\end{coro}

Here $G(z)$ is the Barnes function, which is basically defined by the properties $G(z+1)=G(z)\G(z)$, and $G(0)=G(1)=G(2)=G(3)=1$ (and at all negative integers there are zeros). See Adamchik~\cite{adamchik:contrib-barnes} for a modern exposition and new results about this function.

\subsection{Analytic formula for the volume}

As we did with the supersphere and the complex projective superspace, we can introduce an analytic function expressing the volume.
Define, for complex variables $z$ and $w$,
\begin{equation}
    \V(\sti ; R,z,w):= R^{w(2z-w)}\,(2\sqrt{\pi})^w \frac{G(z-w+1)}{G(z+1)}\,,
\end{equation}
Then, for the normalized volume, we have (with $D=\dim \st_{r}(\C{n|m})= r(2n-r)|mr$)
\begin{equation}
    \frac{\vol \st_{r}(\C{n|m})}{g_D}=     \V(\sti ; R,z,w) \,,
\end{equation}
where    we substitute $z=n-m$, $w=r$\,.

\subsection{The case of the Grassmann supermanifold}

For the complex Grassmann supermanifold $\gr_{r|s}(\C{n|m})=\gra(n|m, r|s)$, we have
\begin{equation*}
    \gr_{r|s}(\C{n|m})\cong \st_{r|s}(\C{n|m})/ \un(r|s)\cong \un(n|m)/\un(n-r|m-s)\times \un(r|s)\,.
\end{equation*}
As in the ordinary case, $\gr_{r|s}(\C{n|m})$ possesses   an analog of the Fubini--Study metric. Dimension and index:
\begin{align*}
    \dim \gr_{r|s}(\C{n|m})&= 2(r|s)\cdot \bigl((n|m)-(r|s)\bigr)\\
    \ind \gr_{r|s}(\C{n|m}) &= 2(r-s) \bigl(n-m-r+s\bigr)\,.
\end{align*}
Define, for complex variables $z$ and $w$,
\begin{equation*}
    \V(\gra ; R,z,w):= R^{2w(z-w)}\,\frac{G(w+1)G(z-w+1)}{G(z+1)}\,,
\end{equation*}
where $G(z)$ is the {Barnes function}.

We expect that the volume of the complex Grassmann supermanifold is given by the following formula:
\begin{equation*}
    \frac{\vol (\gr_{r|s}(\C{n|m}))}{g_{D}}=  \V(\gra ; R,z,w)
\end{equation*}
where $D=\dim \gr_{r|s}(\C{n|m})$, $z=n-m$ and $w=r-s$\,.

\section{Conclusions. Discussion}\label{sec.further}

We considered the volumes of classical supermanifolds: the supersphere, the complex Stiefel and the complex Grassmann supermanifolds. We introduced the notion of normalized volume by isolating the Gaussian factor and showed that the normalized volumes of these supermanifolds are obtained by an analytic continuation from the values for the corresponding ordinary manifolds. We gave a geometric proof for the formula for the volume of the complex Stiefel   supermanifold, which covers Berezin's theorem about the vanishing of the volume of the unitary supergroup as a particular case. On the other hand, our calculations simple examples of compact symplectic supermanifolds with nonzero phase volume (thus disproving Witten's conjecture).

Volumes of supermanifolds exhibit a scale invariance with the scale factor $R^{\ind}$ if the metric is multiplied by $R$. This shows that the volume is scaling invariant for $\ind M=0$. Complex variables in the analytic formulas for the volumes   are also of `index type'.   This, together,   hints at a possibility of extending the theory of volumes to the suitable infinite-dimensional case.

One of the questions that needs to be clarified, is how to agree the analytic expressions when they exist with the special cases such as when the volume identically vanishes for a whole range of parameters. A similar question arises in~\cite{mkrtchyan:nonpertchernsimons}, \cite{mkrtchyanveselov:volumes}, \cite{mkrtchyan:quadrbarnes} in the context of `universal Lie algebra' theory. Elucidating these connections seems   very attractive.

\appendix

\section{Background on supermanifolds, volumes and integrals}\label{sec.prel}

We refer the reader to the following sources concerning supermanifold theory:~ \cite{berezin:antieng},  \cite{leites:bookeng}, \cite{manin:gaugeeng}, \cite{sgelman:eng}, \cite{tv:gitnew}, \cite{deligne:andmorgan}, and \cite{rogers:book}.

\subsection{Tensors on supermanifolds}
One difference with the ordinary case is that it is necessary to distinguish between   `even' and `odd' differentials and, more generally, consider besides the tangent bundle $TM$ and the cotangent bundle $T^*M$   also the `antitangent' bundle  $\Pi TM$ and the `anticotangent' bundle  $\Pi T^*M$ obtained by parity reversion in the fibers.

If the differential $\d f$ of a function  is understood in the most straightforward way,\footnote{\,The differential of $f$ at $x_0$ can be defined, for example,  as  $f-f(x_0) \mod \mathfrak{m}_{x_0}^2$, where $\mathfrak{m}_{x_0}$ denotes the ideal of functions vanishing at $x_0$.}  then the  mapping $\d$ sending  $f$ to   $\d f$ does not change parity, $\widetilde{\d f}=\ft$\,.  (We use here the letter $\d$, not $d$, on purpose.) This is the   \emph{even}  differential of $f$, taking values in $T^*M$. For a tangent vector $X$, $\d \!f (X)= (-1)^{\ft \Xt}\p_Xf$\,.
The differentials $\d x^a$ of  local coordinates $x^a$   make a basis of sections for  $T^*M$. Note    that $\d x^a$ is the \emph{right} dual basis for the basis $\p_a=\lder{}{x^a}$ in $TM$,
\begin{equation*}
    \d x^a (\p_b)=\langle \d x^a, \p_b \rangle= (-1)^{\at}\,\d^a_b
    \quad \text{and} \quad
   \langle \p_a, \d x^b \rangle =(-1)^{\at\bt}\,\langle \d x^b, \p_a \rangle=   \d_a^b\,.
\end{equation*}
The \emph{odd} differential $d$ is defined as the composition $d:=\Pi \circ \d$ and therefore    $dx^a$ make a basis in $\Pi T^*M$. We have $\widetilde{df}=\ft+1$ for all $f$ and $\widetilde{dx^a}=\at+1$. There is an odd pairing $\langle df, X\rangle:=\d \!f (X)$. Hence
\begin{equation*}
    \langle dx^a, \p_b \rangle= \langle \p_b, dx^a \rangle= (-1)^{\at}\,\d^a_b\,.
\end{equation*}

A covariant tensor on a supermanifold is expanded over the tensor products of both $\d x^a$ and $dx^a$. However, products of $dx^a$ may be converted into products of  $\d x^a$, and other way round (by using the sign rule and with a possible appearance of the factor of $\Pi$). Algebraically this corresponds to the  natural isomorphism  $(\Pi V)^{\otimes k}\cong \Pi^k V^{\otimes k}$ for a module $V$. Note also the natural isomorphism $\Pi V^*\cong (\Pi V)^*$.

\begin{example} \label{ex.covtens}
Suppose $T(\p_{a_1},\ldots,\p_{a_k})=T_{a_1\dots a_k}$\,. Then
\begin{equation*}
    T=\d x^{a_1}\otimes \dots \otimes \d x^{a_k} \;T_{a_1\dots a_k} (-1)^{\at_1(\at_2+\ldots+\at_k+\Tt)+\at_2(\at_3+\ldots+\at_k+\Tt)+\ldots+\at_k\Tt}\,.
\end{equation*}
\end{example}

\begin{example} For covariant tensors of rank $2$, we have $dx^a\otimes dx^b=(\Pi \d x^a)\otimes (\Pi \d x^b)=(-1)^{\at} \Pi^2\,\d x^a\otimes \d x^b=(-1)^{\at}  \d x^a\otimes \d x^b$.
\end{example}

The symmetric and the exterior algebra    are defined  as usual (taking into account   the sign rule). The same concerns   symmetrization, alternation, symmetric and wedge products. (Anti)symmetric tensors can be expressed via the symmetric (resp., exterior) products of basis vectors or covectors.

\begin{example} Suppose a covariant tensor as in Example~\ref{ex.covtens} is symmetric. Then
\begin{align*}
    T & = \d x^{a_1}\otimes \dots \otimes \d x^{a_k} \;T_{a_k\dots a_1} (-1)^{(\at_1+ \ldots+\at_k)\Tt}\\
    &=\d x^{a_1}  \dots   \d x^{a_k} \;T_{a_k\dots a_1} (-1)^{(\at_1+ \ldots+\at_k)\Tt}\,,
\end{align*}
where in the second line there is the symmetric product.
\end{example}

For any module $V$,  the natural isomorphism $(\Pi V)^{\otimes k}\cong \Pi^k V^{\otimes k}$ induces a natural isomorphism $S^k(\Pi V)\cong \Pi^k\L^k(V)$. Therefore, a $k$-form  on a supermanifold  can be regarded either as a section of $\L^k(T^*M)$ or as a section of $S^k(\Pi T^*M)$ (with the   parity shift if $k$ is odd).

\begin{example}
Using the identification $dx^a\otimes dx^b= (-1)^{\at}  \d x^a\otimes \d x^b$, we obtain
\begin{multline*}
    dx^a dx^b=\frac{1}{2}\,\left(dx^a\otimes dx^b +(-1)^{(\at+1)(\bt+1)}dx^b\otimes dx^a\right)=\\
    (-1)^{\at}\,\frac{1}{2}\left(\d x^a\otimes \d x^b -(-1)^{\at\bt}\d x^b\otimes \d x^a\right)=
    (-1)^{\at}  \d x^a\wed \d x^b
\end{multline*}
We may write a $2$-form $\o$ as  $\o=dx^adx^b\o_{ba}$, where $\o_{ab}=(-1)^{(\at+1)(\bt+1)}\o_{ba}$, so
\begin{equation*}
    \o=\d x^a\wed \d x^b (-1)^{\at}\o_{ba}\,.
\end{equation*}
Note that for $\o'_{ab}:=(-1)^{\bt}\o_{ab}$ we have $\o'_{ab}=-(-1)^{\at\bt}\o'_{ba}$\,.
\end{example}

The advantage of using $S^k(\Pi T^*M)$ is that $k$-forms in this language can be viewed as fiberwise polynomial functions on the supermanifold $\Pi TM$. Arbitrary functions on $\Pi TM$ are by definition \emph{pseudodifferential forms} on $M$.

\subsection{Berezin integral and volume elements}

Recall  Berezin's change of variables formula.  For an invertible change of variables $x=x(x',\x')$, $\x=\x(x',\x')$ on $\R{n|m}$,
\begin{equation*}
    \int\limits_{\R{n|m}} \!\!\!D(x,\x)\, f(x,\x) =\pm\int\limits_{\R{n|m}} \!\!\!D(x',\x')\,\D{(x,\x)}{(x',\x')}(x',\x')\cdot f\bigl(x(x',\x'),\x(x',\x')\bigr)\,,
\end{equation*}
where $\pm=\sign \det \der{x}{x'}\,(x',0)$,
\begin{equation*}
    \D{(x,\x)}{(x',\x')}= \Ber \der{(x,\x)}{(x',\x')}\,.
\end{equation*}
Here  the function $\Ber$, the  \emph{`superdeterminant'} or \emph{`Berezinian'}, is defined for   even block matrices by the expression
 \begin{equation*}
    \Ber \begin{pmatrix}
           J_{00} & J_{01} \\
           J_{10} & J_{11} \\
         \end{pmatrix}
         :=
         \frac{\det(J_{00}-J_{01}J_{11}^{-1}J_{10})}{\det J_{11}}= \frac{\det J_{00}}{\det(J_{11}-J_{10}J_{00}^{-1}J_{01})}\,.
 \end{equation*}

\begin{Rems} For the change of variables formula to hold, it is essential that all the coefficients of the expansion of $f$ in odd variables rapidly decrease (although only $f_{1\ldots \,m}$ is explicitly present in the definition of the integral).
\end{Rems}

The Berezinian is a multiplicative  function on invertible even matrices and is it essentially uniquely defined by this property. The fact that it is  a rational  expression, not a polynomial, is fundamental. Berezinian enjoys the invariance under the following superanalog of transpose:
\begin{equation*}
    \Ber J^{\,T}=\Ber J
\end{equation*}
for
\begin{equation*}
    \begin{pmatrix}
           J_{00} & J_{01} \\
           J_{10} & J_{11} \\
         \end{pmatrix}^{\!T}=
    \begin{pmatrix}
           J_{00}^t & J_{10}^t \vspace{1pt}\\
           -J_{01}^t & J_{11}^t \\
         \end{pmatrix}\,,
\end{equation*}
where $A^t$ denotes the ordinary transpose of a matrix $A$.

The symbol
\begin{equation*}
    D(x,\x)=D(x^1,\ldots,x^n,\x^1,\ldots,\x^m)
\end{equation*}
is a \emph{coordinate volume element} in the supercase. By definition, this symbol is defined for each coordinate system on $\R{n|m}$ and under a change of coordinates is multiplied by the corresponding Berezinian:
\begin{equation*}
    D(x,\x) = D(x',\x')\,\D{(x,\x)}{(x',\x')}\,.
\end{equation*}
The symbol $D(x,\x)$ is the analog of $d^nx$ for ordinary multiple integrals (over even variables).
In the ordinary case, it is very helpful to  express a coordinate volume element such as $d^nx$ as the complete exterior product $dx^1\wed \ldots \wed dx^n$ of the differentials of coordinates. In the supercase, this is not  possible  directly because Berezinian is   a fraction and cannot arise from  multilinear operations; however, as  a replacement  one can introduce the symbolic bracket
\begin{equation*}
    [dx^1,\ldots,dx^n\,|\,d\x^1,\ldots,d\x^m]\,,
\end{equation*}
as  an alternative notation for $D(x^1,\ldots,x^n,\x^1,\ldots,\x^m)$. 
More generally, for an arbitrary free module $E$ over a commutative superalgebra, we introduce such a `square bracket symbol' as a function of a   basis defined by the following properties:
\begin{itemize}
  \item homogeneity: if a basis element is multiplied by an invertible factor, then the bracket is multiplied by the same factor in the power $+1$ for a basis element in an even position and in the power $-1$ for a basis element in an odd position\footnote{There is an intentional  ambiguity here: we assume that either all even basis elements are   in `even positions' and all odd basis elements are in `odd positions', or the other way round. This allows us to suppress the difference between $dx^a$, $d\x^{\mu}$ and $\d x^a$, $\d \x^{\mu}$ in the notation for the bracket.};
  \item invariance under elementary transformations: the symbol does not change when a basis element is replaced by the sum with another   element with a coefficient of the appropriate parity.
\end{itemize}
These properties model the characteristic properties of Berezinian. The `square bracket' of a basis of a free module $E$ gives  a basis element of the one-dimensional module  $\Ber E$ (the Berezinian of a free module).
One can  learn quickly that it is as convenient to make calculations with the symbol $[dx^1,\ldots,dx^n\,|\, d\xi^1,\ldots,d\xi^m]$ on supermanifolds as with the exterior product $dx^1\wed \ldots \wed dx^n$ on ordinary manifolds. In the following we may use self-explanatory abbreviated notations such as $[dx|d\x]$ for $[dx^1,\ldots,dx^n\,|\, d\xi^1,\ldots,d\xi^m]$.

\begin{example} Consider $\C{n|m}$ with complex coordinates $z^a=x^a+iy^a$, $\z^{\mu}=\x^{\mu}+i\h^{\mu}$. By applying elementary transformations, we obtain $[dz,\bar{dz}\,|\, d\z,\bar{d\z}]= [dx+idy,dx-idy \,|\, d\x+id\h,d\x-id\h]=[2dx,dx-idy \,|\, 2d\x,d\x-id\h]=[2dx,-idy \,|\, 2d\x,-id\h]$. By homogeneity, we arrive at
\begin{equation*}
    [dz,d\bar{z}\,|\, d\z,d\bar{\z}]= (-2i)^{n-m}[dx,dy\,|\,d\x,d\h]\,,
\end{equation*}
hence
\begin{equation*}
    [dx,dy\,|\,d\x,d\h]=\left(\frac{i}{2}\right)^{\!n-m}[dz,d\bar{z}\,|\, d\z,d\bar{\z}]\,.
\end{equation*}
\end{example}

\subsection{Gaussian integrals. Gaussian factor}

Consider a non-degenerate even quadratic function $Q(x)$ on $\R{n|2m}$. That the odd part of the dimension has to be $2m$ is required by the non-degeneracy. It corresponds to an even symmetric bilinear form, which we denote by the same letter, so that $Q(x)=Q(x,x)$. In coordinates,
\begin{equation*}
    Q(x)=x^ax^bQ_{ba}=x^aQ_{ab}(-1)^{\bt}x^b\,.
\end{equation*}
where $Q_{ab}=Q_{ba}(-1)^{\at\bt}$.
Consider the Gaussian integral defined by $Q$. One can see that
\begin{equation*}
    \int\limits_{\R{n|2m}} e^{-Q(x)} Dx =\frac{1}{\left(\Ber Q\right)^{1/2}}\, (\sqrt{\pi})^n2^m\,,
\end{equation*}
where $\Ber Q$ is the Berezinian of the matrix $(Q_{ab})$. Here $(\sqrt{\pi})^n2^m$ is the value of the `standard' Gaussian integral corresponding to the matrix
\begin{equation*}
    \begin{pmatrix}
      I_n & 0 \\
      0 & J_{2m} \\
    \end{pmatrix}\,,
\end{equation*}
with $I_n$ being the identity matrix and $J_{2m}=\diag\left(\left(\begin{smallmatrix}0 & -1\\ 1 & 0\end{smallmatrix}\right), \ldots, \left(\begin{smallmatrix}0 & -1\\ 1 & 0\end{smallmatrix}\right)\right)$\,,
for which
\begin{equation*}
    Q(x)=(x^1)^2+\ldots+(x^n)^2 + 2\x^1\x^2+\ldots + 2\x^{2m-1}\x^{2m}
\end{equation*}
(here we have switched to separate notations for even and odd variables).

We number $(\sqrt{\pi})^n2^m$, which is a function of dimension, will play a role in our statements. We shall denote it
\begin{equation*}
    g_{n|2m}:=(\sqrt{\pi})^n2^m\,,
\end{equation*}
and refer to as the \emph{Gaussian factor}. (Formally, we can use the same notation for an arbitrary $n|m$.)

\section{Volumes and Riemannian submersions}

\subsection{Volume elements arising from Riemannian and symplectic structures} \label{subsec.volel} On a supermanifold,  an even Riemannian metric or  even symplectic form induce volume elements in the same way as   Riemannian and symplectic structures for ordinary manifolds.\footnote{\,`Even' or `odd' here refers to the parity of the corresponding bilinear form. Odd structures behave quite differently from their classical prototypes. In particular, there are no volume elements associated with them, which has important consequences. Odd symplectic geometry underlies the Batalin--Vilkovisky quantization, while odd Riemannian geometry has so far attracted little attention. See more, for example, in~\cite{tv:laplace1}.}

Consider an even covariant tensor $T$ of rank $2$ (without any conditions of symmetry). In local coordinates,
\begin{equation*}
   T=  \d x^a \otimes \d x^b \,T_{ab}(-1)^{\at\bt}\,,
\end{equation*}
where $T_{ab}=T(e_a,e_b)$, and $\widetilde{T_{ab}}=\at+\bt$. Under a change of coordinates,
\begin{align*}
    T_{a'b'} &= \der{x^{a}}{x^{a'}}\,\der{x^{b}}{x^{b'}}\,(-1)^{\at(\bt+\bt')}\,T_{ab} \\
    &=\der{x^{a}}{x^{a'}}\,T_{ab}\,\der{x^{b}}{x^{b'}}\,(-1)^{\bt(\bt'+1)}
\end{align*}
(where we recognize the supertranspose of the Jacobi matrix).
It follows that
\begin{equation*}
    \Ber (T_{a'b'})= \Ber (T_{ab})\cdot \left(\Ber \!\left(\der{x^{a}}{x^{a'}}\right)\right)^{\!2}\,.
\end{equation*}
Hence
\begin{equation*}
    dV:=\sqrt{\Ber (T_{ab})}\,Dx
\end{equation*}
is an invariant volume element, regardless of the symmetry of $T$.

Suppose an even Riemannian metric is given:
\begin{equation*}
    \d s^2 = \d x^a \d x^b \,g_{ba}= \d x^a g_{ab}\, (-1)^{\bt}\d x^b
\end{equation*}

(even) on a supermanifold has an expression
\begin{equation*}
    \d s^2 = \d x^a \d x^b \,g_{ba}= \d x^a g_{ab}\, (-1)^{\bt}\d x^b
\end{equation*}
in local coordinates. Here $g_{ab}=g_{ab}(x)$, $\widetilde{g_{ab}}=\at+\bt$, and $g_{ab}=(-1)^{\at\bt}g_{ba}$. We use $\d$ to denote the even differential (so that, in particular, $\widetilde{\d x^a}=\at$). Under a change of coordinates,
\begin{equation*}
    g_{a'b'} = (-1)^{\at(\bt+\bt')}\der{x^{a}}{x^{a'}}\,\der{x^{b}}{x^{b'}}\,g_{ab}=
    \der{x^{a}}{x^{a'}}\,g_{ab}\,(-1)^{\bt(\bt'+1)}\der{x^{b}}{x^{b'}}\,.
\end{equation*}
Hence, for $g:=\Ber (g_{ab})$, we obtain
\begin{equation*}
    g' = g\cdot \left(\Ber \!\left(\der{x^{a}}{x^{a'}}\right)\right)^{\!2}\,,
\end{equation*}
where $g'=\Ber (g_{a'b'})$, and therefore
\begin{equation*}
    dV:=\sqrt{g}\,Dx
\end{equation*}
is an invariant volume element in  complete similarity with the ordinary case.

Consider now an even symplectic structure. One can similarly obtain   an induced volume element with the same expression as in the Riemannian case. It is sometimes referred to it as   the `Liouville volume element'.

\begin{Rem} For   ordinary manifolds, the Liouville volume element can be expressed as the top-degree form $\o^n/n!$\,, where $\dim M=2n$. The analog of this expression in the supercase is the pseudodifferential form $e^{-\o}$ (see, e.g.~\cite{tv:gitnew}).  When we integrate such objects, we have to integrate over   $x^a$ and $dx^a$ (treated as independent variables). In particular, integrating out the variables $dx^a$ in $e^{-\o}$ gives the square root of the Berezinian of $\o_{ab}$. The difference with the ordinary case is that this square root is no longer a polynomial in   $\o_{ab}$.
\end{Rem}

It is known that the Riemannian and symplectic structures come together in the notion of a K\"{a}hler structure. In the supercase, we have a tensor $H$ defining an even Hermitian metric on a complex supermanifold with holomorphic local coordinates $z^a$,

\begin{equation*}
    H=\d z^a \otimes \d \bar{z^b} \,h_{a\bar b}(-1)^{\at\bt}
\end{equation*}
where $H_{a\bar b}=(-1)^{\at\bt}\widebar{H_{b\bar a}}$. Then $\d s^2:=\Re H$ gives a Riemannian structure and $\o:=-\frac{1}{2} \Im H$ gives a symplectic structure. The corresponding volume elements (Riemannian and symplectic) coincide.

\subsection{Riemannian submersions}

Recall the following classical notion. Consider a fiber bundle $p\co E\to M$ where both $E$ and $M$ are endowed with   Riemannian metrics. The tangent bundle $TE$ can be decomposed into the direct sum $VE\oplus HE$ where  $VE=\Ker T\!p$  (the vertical subbundle)  and the horizontal subbundle $HE=(VE)^{\perp}$   is defined as the orthogonal complement w.r.t. the metric on $E$. The restriction of the tangent map $T\!p$ maps each horizontal subspace $HE_z$  isomorphically onto  the tangent space $T_{p(z)}M$.   The fiber bundle  is called a \emph{Riemannian submersion} if   $T\!p(z)_{|HE_z}\co HE_z \to T_{p(z)}M$ is an isometry (for all $z\in E$).
This notion readily extends to supermanifolds.

Consider the direct product coordinates $x^a,y^i$ on $E$ so that $x^a$ are local coordinates on $M$ and the bundle projection maps $(x^a,y^i)$ to $(x^a)$. The Riemannian metrics  on $E$ and $M$  can be written as
\begin{equation*}
    \d s^2_E=\d x^a \d x^b\,g_{ba}(x,y) + 2 \d x^a \d y^i \,g_{ia}(x,y)+ \d y^i\d y^j \,g_{ji}(x,y)
\end{equation*}
and
\begin{equation*}
    \d s^2_M=\d x^a \d x^b\, g_{ba}^0(x)\,.
\end{equation*}

\begin{prop} The fiber bundle $p\co E\to M$ is a Riemannian submersion if there is the following relation for the components of the metrics
\begin{equation*}
      g_{ab}(x,y)=g_{ab}^0(x) +g_{ak}(x,y)g^{kl}(x,y)g_{lb}(x,y)
\end{equation*}
\end{prop}

\begin{cor} For the Berezinian  of the metric tensor,
\begin{equation*}
    \Ber\begin{pmatrix}
          g_{ab} & g_{aj} \\
          g_{aj} & g_{ij} \\
        \end{pmatrix}=
       \Ber\begin{pmatrix}
          g_{ab}^0 & 0 \\
          g_{aj} & g_{ij} \\
        \end{pmatrix}\,.
\end{equation*}
\end{cor}
This implies the factorization of the volume element of $E$ at a given point in the product of the volume element of $M$ (at the corresponding point) and the volume element of the fiber through this point.

Therefore there is a `Cavaliery principle' for Riemannian submersions: the volume of the total space is the integral of the volumes of the fibers over the base. If all fibers have the same volume (such as in a homogeneous situation), the volume of the total space is the product of the volume of the base and the volume of the fiber.



\def\cprime{$'$} \def\cprime{$'$} \def\cprime{$'$} \def\cprime{$'$}
  \def\cprime{$'$} \def\cprime{$'$} \def\cprime{$'$} \def\cprime{$'$}
  \def\cprime{$'$} \def\cprime{$'$} \def\cprime{$'$} \def\cprime{$'$}
  \def\cprime{$'$}

\end{document}